\documentclass[a4paper,11pt,reqno]{amsart}

\usepackage{spiros}
\usepackage[pdftex]{graphicx}
\usepackage[all]{xy}

\newtheorem{thm1}{Theorem}
\newtheorem{prop1}[thm1]{Proposition}

\title{A controlled local-global theorem for simplicial complexes}
\author{Spiros Adams-Florou}
\newtheorem*{ack*}{Acknowledgement}

\begin{document}

\maketitle

\begin{abstract}
In this paper we prove that a simplicial map of finite-dimensional locally finite simplicial complexes has contractible point inverses if and only if it is an $\ep$-controlled homotopy equivalence for all $\ep>0$ if and only if $f\times \id_\R$ is a bounded homotopy equivalence measured in the open cone over the target. This confirms for such a space $X$ the slogan that arbitrarily fine control over $X$ corresponds to bounded control over the open cone $O(X_+)$. For the proof a one parameter family of cellulations $\{X_\ep^\prime\}_{0<\ep<\ep(X)}$ is constructed which provides a retracting map for $X$ which can be used to compensate for sufficiently small control.
\end{abstract}

%QEDs, A note on properness of maps and simplicial approximations

\section{Introduction}
%Currently poor intro, more history, relevance, etc, explain what simple means, talk about approximation theorems, talk about simple htpy equivs.

% The idea of using estimates in geometric topology dates back to Connell and Hollingsworth's paper \cite{ConnHoll} in which the authors introduce estimates in order to compute algebraic obstructions. After this, controlled topology was developed by Chapman (\cite{Chapcontrol}), Ferry (\cite{Ferryep}) and Quinn (\cite{Quinn1}, \cite{Quinn2}). The general approach of controlled topology is to put an estimate on a geometric obstruction (usually by use of a metric). It is often possible to prove that if the size of the estimate is sufficiently small then the obstruction must vanish.

A homeomorphism has point inverses which are all points. If a map $f$ is homotopic to a homeomorphism it is reasonable to suppose that $f$ might have point inverses that are `close' to being points in some suitable sense. Controlled topology takes `close' to mean small with respect to a metric. One then studies maps with small point inverses and attempts to prove that such a map is homotopic to a homeomorphism. 

This approach has many successes in the literature: as a consequence of Chapman and Ferry's $\alpha$-approximation theorem (\cite{chapfer}) a map between closed metric topological manifolds with sufficiently small point inverses is homotopic to a homeomorphism through maps with small point inverses. One can also consider maps where the point inverses all have the homotopy groups of a point, i.e. are contractible. In the non-manifold case Cohen proves in \cite{cohen} that a p.l. map of finite polyhedra with contractible point inverses is a simple homotopy equivalence. 

When doing controlled topology it is desirable that the space we consider, $X$, comes equipped with a metric. In the absence of a metric it is sufficient that $X$ has at least a map $p:X\to M$ to a metric space $(M,d)$, called a \textit{control map}, which then allows us to measure distances in $M$. In general to be able to detect information about $X$ from the control map and the metric on $M$ we would ideally like $p$ to be highly connected.

Let $f:X\to Y$ be a map of spaces equipped with control maps $p:X\to M$, $q:Y\to M$ to a metric space $(M,d)$. We say that $f:(X,p) \to (Y,q)$ is \textit{$\ep$-controlled} if $f$ commutes with the control maps $p$ and $q$ up to a discrepancy of $\ep$, i.e. for all $x\in X$, $d(p(x), qf(x))<\ep$. We say that $f:(X,p) \to (Y,q)$ is an \textit{$\ep$-controlled homotopy equivalence} if there exists a homotopy inverse $g$ and homotopies $h_1:g\circ f \sim \id_X$ and $h_2:f\circ g\sim \id_Y$ such that all of $f:(X,p)\to (Y,q)$, $g:(Y,q)\to (X,p)$, $h_1:(X\times\R,p\times\id_\R) \to (X,p)$ and $h_2:(Y\times\R,q\times\id_\R)\to (Y,q)$ are $\ep$-controlled maps. %As before we call $\ep$ the \textit{control} of the homotopy equivalence $f$. 

Note that an $\ep$-controlled homotopy equivalence $f$ and its inverse $g$ do not move points more than a distance $\ep$ when measured in $M$ and that the homotopy tracks are no longer than $\ep$ when measured in $M$. If $X$ and $Y$ are also metric spaces it is perfectly possible that the homotopy tracks are large in $X$ or $Y$ and only become small after mapping to $M$.

Controlled topology is not functorial because the composition of two maps with control less than $\ep$ is a map with control less than $2\ep$. This motivates Pedersen's development in \cite{PedKminusi} and \cite{PedInvariants} of bounded topology, where the emphasis is no longer on how small the control is but rather just that it is finite. A map $f:(X,p) \to (Y,q)$ is called \textit{bounded} if $f$ commutes with the control maps up to a finite discrepancy $B$. Similarly a bounded homotopy equivalence is one where all the maps and homotopies are bounded. Bounded topology is functorial as the sum of two finite bounds remains finite. 

%Typically controlled topology is concerned with small control measured in a compact space, whereas bounded topology is concerned with finite bounds measured in a non-compact space. Bounded topology is coarse in the sense that local topology is invisible to the control map; only the global behaviour out towards infinity is noticed. Controlled topology detects things with non-zero size, but is blind to the arbitrarily small, for example a missing point in $X$ cannot be detected. 
In \cite{epsurgthy} Ferry and Pedersen suggest a relationship between controlled topology on a space $X$ and bounded topology on the open cone $O(X_+)$ when they write in a footnote:
\begin{quotation}
``It is easy to see that if $Z$ is a Poincar\'{e} duality space with a map $Z \to K$ such that $Z$ has $\ep$-Poincar\'{e} duality for all $\ep>0$ when measured in $K$ (after subdivision), e.g. a homology manifold, then $Z\times\R$ is an $O(K_+)$-bounded Poincar\'{e} complex. The converse (while true) will not concern us here.'' 
\end{quotation}
% Morally, since the metric on $O(K_+)$ scales by a factor of $t$ at height $t$, a finite bound $B$ on the open cone corresponds to a bound proportional to $\frac{B}{t}$ on $(M\times\{t\},d)$.
For $X$ a proper subset of $S^n$ the open cone $O(X_+)\subset \R^{n+1}$ is the union of all rays from the origin $0\in \R^{n+1}$ through points in $X_+=X\sqcup \{x_0\}$ together with the subspace metric. There is a natural map $j_X:X\times\R \to O(X_+)$, called the \textit{coning map}, given by \[j_X(x,t) := \brcc{tx,}{t\geqslant 0,}{-tx_0,}{t\leqslant 0.}\]See section \ref{prelims} for a more general definition of the open cone and the coning map for more general metric spaces. 

The footnote above leads one to conjecture that $f:(X,qf) \to (Y,q)$ is an $\ep$-controlled homotopy equivalence for all $\ep>0$ if and only if \[f\times \id_\R:(X\times\R, j_Y(qf\times \id_\R)) \to (Y\times\R,j_Y)\] is a bounded homotopy equivalence. In this paper we prove this conjecture for the case of a simplicial map of finite-dimensional locally finite (henceforth f.d. l.f.) simplicial complexes measured in the target. We may measure in the target since such complexes come naturally equipped with a path metric. We prove

\begin{thm1}\label{maintopthm}
Let $f:X\to Y$ be a simplicial map of f.d. l.f. simplicial complexes with $Y$ equipped with the path metric. Then the following are equivalent:
\begin{enumerate}[(i)]
 \item $f$ has contractible point inverses,
 \item $f:(X,f)\to (Y,\id_Y)$ is an $\ep$-controlled homotopy equivalence for all $\ep>0$,
 \item $f\times\id_\R: (X\times\R, j_Y(f\times\id_\R)) \to (Y\times\R,j_Y)$ is a bounded homotopy equivalence.
\end{enumerate}
\end{thm1}

%This result can be considered an extension of GORUAWOGUR
Working with simplicial maps makes life much easier - one needs only check that the point inverses of the barycentres are contractible:

\begin{prop1}\label{two}
Let $f:X\to Y$ be a simplicial map of l.f. f.d. simplicial complexes. Then 
\begin{enumerate}[(i)]
 \item for all simplices $\sigma \in Y$, there is a p.l. isomorphism $f^{-1}(\mathring{\sigma}) \cong f^{-1}(\widehat{\sigma}) \times \mathring{\sigma},$
 \item $f$ has contractible point inverses if and only if $f^{-1}(\widehat{\sigma})$ is contractible for all $\sigma\in Y$.
\end{enumerate}
\end{prop1}

Moreover, simplicial maps allow us to `lift' certain properties of the target space to the preimage, in particular the fact that open stars deformation retract onto open simplices:

\begin{prop1}\label{three}
Let $f:X\to Y$ be a simplicial map of f.d. l.f. simplicial complexes. Then for all $\sigma\in Y$, $f^{-1}(\mathrm{st}(\sigma))$ p.l. deformation retracts onto $f^{-1}(\mathring{\sigma})$.\footnote{By a deformation retract we mean a strong deformation retract.}
\end{prop1}

If, as in the theorem, we additionally suppose that a simplicial map $f:X\to Y$ has contractible point inverses, then $f$ turns out to have the \textit{approximate homotopy lifting property}: for all $\ep>0$, the lifting problem

\begin{displaymath}
 \xymatrix{
  Z\times \{0\} \ar[r]^-{h} \ar[d] & X\ar[d]^-{f} \\
  Z\times I \ar@{-->}[ur]^-{\widetilde{H}_\ep} \ar[r]^-{H} & Y
 }
\end{displaymath}
has a solution $\widetilde{H}_\ep: Z\times I \to X$ such that the diagram commutes up to $\ep$, i.e. \[d_Y(H(z,t),f(\widetilde{H}_\ep(z,t))<\ep\] for all $(z,t)\in Z\times I$, where $d_Y$ is the metric on $Y$. This is precisely the definition of an approximate fibration given by Coram and Duvall in \cite{CoramDuvall}.

The key ingredient in proving $(i)\Rightarrow(iii)$ and obtaining the approximate homotopy lifting property is the construction and use of the \textit{fundamental $\ep$-subdivision cellulation} $X^\prime_\ep$ of an f.d. l.f. simplicial complex $X$. The $X_\ep^\prime$ are a family of cellulations with $\molim_{\ep\to 0}X_\ep^\prime = X$ similar to the family of cellulations obtained by taking slices $||X||\times\{t\}$ of the prism $||X||\times [0,1]$ triangulated so that $||X||\times\{0\}$ is given the triangulation $X$ and $||X||\times\{1\}$ the barycentric subdivision $Sd\, X$. The key difference is that the cellulations $X_\ep^\prime$ are defined in such a way as to guarantee that the homotopy from $X_\ep^\prime$ to $X$ through $X_\delta^\prime$ for $\delta\in (0,\ep)$ has control $\ep$. These cellulations provide retracting maps that compensate for $\ep$-control when proving squeezing results. This is precisely what is missing when trying to prove such results for a more general class of spaces. 

%Rmk that open star is intersection of open stars of the vertices

Section \ref{prelims} recaps some necessary preliminaries. In section \ref{epsubdivcell} the fundamental $\ep$-subdivision cellulation of an f.d. l.f. simplicial complex is defined and a few useful properties explained. In section \ref{proof} Propositions \ref{two} and \ref{three} are proved and consequently a direct proof of Theorem \ref{maintopthm} is given.
\begin{ack*}
This work is partially supported by Prof. Michael Weiss' Humboldt Professorship.
\end{ack*}

\section{Preliminaries}\label{prelims}
In this paper only locally finite finite-dimensional simplicial complexes will be considered. Such a space $X$ shall be given a metric $d_X$, called the \textit{standard metric}, as follows. First define the \textit{standard $n$-simplex} $\Delta^n$ in $\R^{n+1}$ as the join of the points $e_0 = (1,0,\ldots, 0)$, \ldots, $e_n = (0,\ldots,0,1)\in \R^{n+1}$. $\Delta^n$ is given the subspace metric $d_{\Delta^n}$ of the standard $\ell_2$-metric on $\R^{n+1}$. The locally finite finite-dimensional simplicial complex $X$ is then given the path metric whose restriction to each $n$-simplex is $d_{\Delta^n}$. Distances between points in different connected components are thus $\infty$. See $\S 4$ of \cite{bartelssqueezing} or Definition $3.1$ of \cite{HR95} for more details.

Let $p:Y\to X$ be a simplicial map of locally-finite simplicial complexes equipped with standard metrics. For $\sigma$ a simplex in $Y$, the \textit{diameter of $\sigma$ measured in $X$} is \[\mathrm{diam}(\sigma):= \sup_{x,y\in\sigma}{d_X(p(x),p(y))}.\] The \textit{radius of $\sigma$ measured in $X$} is \[\mathrm{rad}(\sigma) := \inf_{x\in\partial\sigma}d_X(p(\widehat{\sigma}),p(x)).\] The \textit{mesh of $X$ measured in $Y$} is \[\mesh(X):= \sup_{\sigma\in X}\{\mathrm{diam}(\sigma)\}.\] The \textit{comesh of $X$ measured in $Y$} is \[\comesh(X):= \inf_{\sigma\in X, |\sigma|\neq 0}\{\mathrm{rad}(\sigma)\}.\]

Using the standard metric on $X$ and $\id_X:X\to X$ as the control map $\diam (\sigma) = \sqrt{2}$ and $\rad (\sigma) = \frac{1}{\sqrt{|\sigma|(|\sigma|+1)}},$ for all $\sigma\in X$, so consequently $\mesh (X) = \sqrt{2}$ and if $X$ is $n$-dimensional $\comesh (X) = \frac{1}{\sqrt{n(n+1)}}$.

The open star st$(\sigma)$ of a simplex $\sigma \in X$ is defined by \[\mathrm{st}(\sigma):= \bigcup_{\tau\geqslant \sigma}\mathring{\tau}.\]

% In this paper we work with finite-dimensional locally finite (henceforth f.d. l.f.) simplicial complexes. Such a space naturally comes with a complete metric: the path metric obtained by given each simplex the subspace metric from its standard embedding in Euclidean space. See $\S 4$ of \cite{bartelssqueezing} or Definition $3.1$ of \cite{HR95} for more details. We also consider the open cone $O(X_+)$ of a f.d. l.f. simplicial complex $X$. 

The open cone was first considered by Pedersen and Weibel in \cite{kthyhom} where it was defined for subsets of $S^n$. This definition was extended to more general spaces by Anderson and Munkholm in \cite{AndMunk}. We make the following definition: For a complete metric space $(M,d)$ the \textit{open cone} $O(M_+)$ is defined to be the identification space $M\times\R/\sim$ with $(m,t)\sim(m^\prime,t)$ for all $m,m^\prime\in M$ if $t\leqslant 0$. We define a metric $d_{O(M_+)}$ on $O(M_+)$ by setting 
\begin{eqnarray*}
 d_{O(M_+)}((m,t),(m^\prime,t)) &=& \brcc{td(m,m^\prime),}{t\geqslant 0,}{0,}{t\leqslant 0,} \\
 d_{O(M_+)}((m,t),(m,s)) &=& |t-s|
\end{eqnarray*}
 and defining $d_{O(M_+)}((m,t),(m^\prime,s))$ to be the infimum over all paths from $(m,t)$ to $(m^\prime,s)$, which are piecewise geodesics in either $M\times\{r\}$ or $\{n\}\times\R$, of the length of the path. I.e.
\[d_{O(M_+)}((m,t),(m^\prime,s)) = \max\{\min\{t,s\},0\}d_X(m,m^\prime) + |t-s|.\]
This metric is carefully chosen so that 
\[ d_{O(M_+)}|_{M\times\{t\}} = \brcc{td_{O(M_+)}|_{M\times\{1\}},}{t\geqslant 0,}{0,}{t\leqslant 0.}\]
This is precisely the metric used by Anderson and Munkholm in \cite{AndMunk} and also by Siebenmann and Sullivan in \cite{SiebSull}, but there is a notable distinction: we do not necessarily require that our metric space $(M,d)$ has a finite bound. 

There is a natural map $j_X:X\times\R \to O(X_+)$ given by the quotient map
\begin{eqnarray*}
X\times\R &\to& X\times\R/\sim \\ 
(x,t) &\mapsto& [(x,t)].
\end{eqnarray*}
We call this the \textit{coning map}.

For $M$ a proper subset of $S^n$ with the subspace metric, the open cone $O(M_+)$ can be thought of as the subset of $\R^{n+1}$ consisting of all the points in the rays out of the origin through points in $M_+:= M \cup \{pt\}$ with the subspace metric. This is not the same as the metric we just defined above but it is Lipschitz equivalent. 

\section{Subdivision cellulations}\label{epsubdivcell}
In this section we construct a controlled $1$-parameter family of subdivision cellulations of $X$ which shall be used later in constructing controlled homotopies. This $1$-parameter family is defined in analogy to the $1$-parameter family of subdivision cellulations obtained by restricting a triangulation of the prism $X\times I$ to the slices $\{ X\times \{t\} \}_{0<t<1}$.

Given an f.d. l.f. simplicial complex $X$ and its barycentric subdivision $Sd\, X$, we may triangulate the prism $||X||\times I$ so that $||X||\times \{0\}$ has triangulation $X$ and $||X||\times\{1\}$ has triangulation $Sd\, X$. % (see \cite{hatcher} page 112) 
\begin{defn}
The \textit{canonical triangulation of $||X||\times I$ from $X$ to $Sd\, X$} is defined to have one $(|\sigma|+n+1)$-simplex \[ (\sigma\times\{0\})*(\widehat{\sigma}_0\ldots \widehat{\sigma}_n\times \{1\})\] in $||X||\times I$ for every chain of inclusions in $X$ of the form \[\sigma\leqslant \sigma_0<\ldots<\sigma_n.\] With a slight abuse of terminology we shall call such a chain of inclusions a \textit{flag in $X$ of length $n$}. It may easily be verified that this indeed gives a triangulation.
\end{defn}

\begin{ex}
Let $X$ be a $2$-simplex. \Figref{Fig:PrismSlicing} illustrates the canonical triangulation of $||X||\times I$ and what the induced cellulations of the slice $||X||\times\{0.5\}$ is.
\begin{figure}[h!]
\begin{center}
{
\psfrag{t=0}[]{$t=0$}
\psfrag{t=0.5}[]{$t=0.5$}
\psfrag{t=1}[]{$t=1$}
\includegraphics[width=9cm]{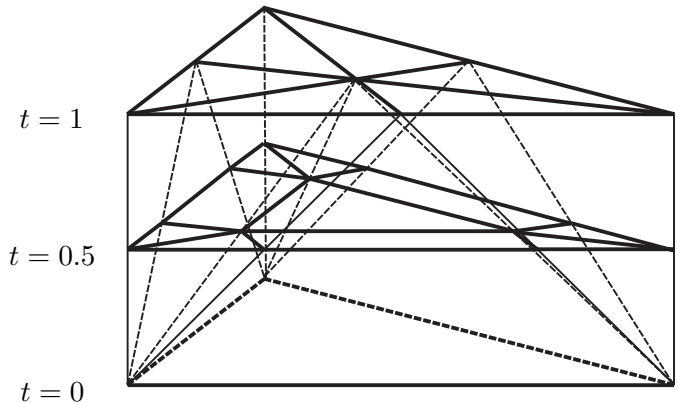}
}
\caption{Obtaining cellulations from the prism.}
\label{Fig:PrismSlicing}
\end{center}
\end{figure}
\end{ex}

The slices $\{ ||X||\times \{t\} \}_{0<t<1}$ form a continuous family of cellulations of $||X||$ from $X$ to $Sd\, X$. Mapping cells identically to corresponding cells and taking the limit as $t\to 0$ there is a straight line homotopy on $||X||$ sending the cellulation of $||X||\times\{t\}$ to $X$ by mapping through the cellulations $(||X||\times \{s\})_{0<s<t}$. We adapt this procedure to give a family of cellulations, $X_\ep^\prime$, where the straight line homotopy from $X_\ep^\prime$ to $X$ has control at most $\ep$ measured in $X$. 

\begin{defn}\label{Defn:flagcellulation}
Define the \textit{flag cellulation of $X$} by \[\chi(X):= \bigcup_{m=0}^{\mathrm{dim}(X)}\bigcup_{\sigma\leqslant\sigma_0 < \ldots < \sigma_m}{\sigma \times \widehat{\sigma}_0\ldots \widehat{\sigma}_m} \subset X\times Sd\, X.\]
\end{defn}
Observe that $\chi(X)$ has the same cellulation as that inherited by $||X||\times \{t\}$ for any $t\in (0,1)$ from the canonical triangulation of the prism from $X$ to $Sd\, X$. We now construct a $1$-parameter family of p.l. isomorphisms $\Gamma_\ep: \chi(X) \to X$ which shall be used to give $X$ a $1$-parameter family of cellulations.  %so that $\Gamma_\ep(v\times D(v,X)) \subset B_\ep(v)$

\begin{defn}\label{Defn:fundamentalcellulation}
For $0\leqslant\ep<\comesh(X)$ define a map $\Gamma_\ep : \chi(X) \to X$ by
\begin{eqnarray*}
 \Gamma_\ep(v\times \widehat{v}) &:=&v, \quad \;\mathrm{for}\;\mathrm{all}\;\mathrm{vertice}\;v\in X, \\
  \Gamma_\ep(v\times \widehat{\tau}) &:=& \partial B_\ep(v)\cap \widehat{v}\widehat{\tau}, \quad\mathrm{for}\;\mathrm{all}\;\mathrm{inclusions}\;\mathrm{of}\;\mathrm{a}\;\mathrm{vertex}\; v< \tau,
\end{eqnarray*}
where $\partial B_\ep(v)$ is the sphere of radius $\ep$ centred at the vertex $v$ and $\partial B_0(v):= v$. 

Extend $\Gamma_\ep$ piecewise linearly over each cell of $\chi$ by
\begin{eqnarray*}
 \Gamma_\ep: \sigma \times \widehat{\sigma}_0\ldots \widehat{\sigma}_m &\to& X \\
 (s_0,\ldots, s_n,t_0,\ldots, t_m) &\mapsto& \sum_{i=0}^{n}\sum_{j=0}^{m}{s_it_j\Gamma_\ep(v_i\times\widehat{\sigma}_j)},
\end{eqnarray*}
where $\sigma = v_0\ldots v_n$ with barycentric coordinates $(s_0,\ldots, s_n)$ and $\widehat{\sigma}_0\ldots \widehat{\sigma}_m$ has barycentric coordinates $(t_0,\ldots, t_m)$. 

We call the image under $\Gamma_\ep$ of the flag cellulation the \textit{fundamental $\ep$-subdivision cellulation of $X$} and denote it by $X_\ep^\prime$. We use the following notation for the cells of $X_\ep^\prime$:
\begin{eqnarray*}
\Gamma_{\sigma_0,\ldots, \sigma_m}(\sigma)&:=& \Gamma_\ep(\sigma\times \widehat{\sigma}_0\ldots \widehat{\sigma}_m), \\
\Gamma_{\sigma_0,\ldots, \sigma_m}(\mathring{\sigma})&:=& \Gamma_\ep(\mathring{\sigma}\times \widehat{\sigma}_0\ldots \widehat{\sigma}_m),
\end{eqnarray*}
for all flags $\sigma\leqslant \sigma_0<\ldots< \sigma_n$.
\qed\end{defn}

\begin{ex}\label{starstar}
Let $X$ be the simplex $\sigma=v_0v_1v_2$ with faces labelled $\tau_0=v_0v_1$, $\tau_1=v_1v_2$ and $\tau_2=v_0v_2$, then the fundamental $\ep$-subdivision cellulation of $X$ is as in \Figref{Fig:fundamentalcellulation}. Each $\Gamma_{\sigma_0,\ldots,\sigma_i}(\tau)$ is the closed cell pointed to by the arrow.
\begin{figure}[tbh!]
\begin{center}
{
\psfrag{v0}{$\rho_0$}
\psfrag{v1}{$\rho_1$}
\psfrag{v2}{$\rho_2$}
\psfrag{ep}{$\ep$}
\psfrag{t1st1}{$\Gamma_{\tau_0,\sigma}(\tau_0)$}
\psfrag{t1t1}{$\Gamma_{\tau_0}(\tau_0)$}
\psfrag{st1}{$\Gamma_{\sigma}(\tau_0)$}
\psfrag{r2t2sr2}{$\Gamma_{\rho_1,\tau_1,\sigma}(\rho_1)$}
\psfrag{r2sr2}{$\Gamma_{\rho_1,\sigma}(\rho_1)$}
\psfrag{t2sr2}{$\Gamma_{\tau_1,\sigma}(\rho_1)$}
\psfrag{r2t2r2}{$\Gamma_{\rho_1,\tau_1}(\rho_1)$}
\psfrag{r3sr3}[tr][tr]{$\Gamma_{\rho_2,\sigma}(\rho_2)$}
\psfrag{sr3}{$\Gamma_{\sigma}(\rho_2)$}
\psfrag{r3r3}{$\Gamma_{\rho_2}(\rho_2)$}
\psfrag{s}[][]{$\Gamma_{\sigma}(\sigma)$}
\includegraphics[width=10cm]{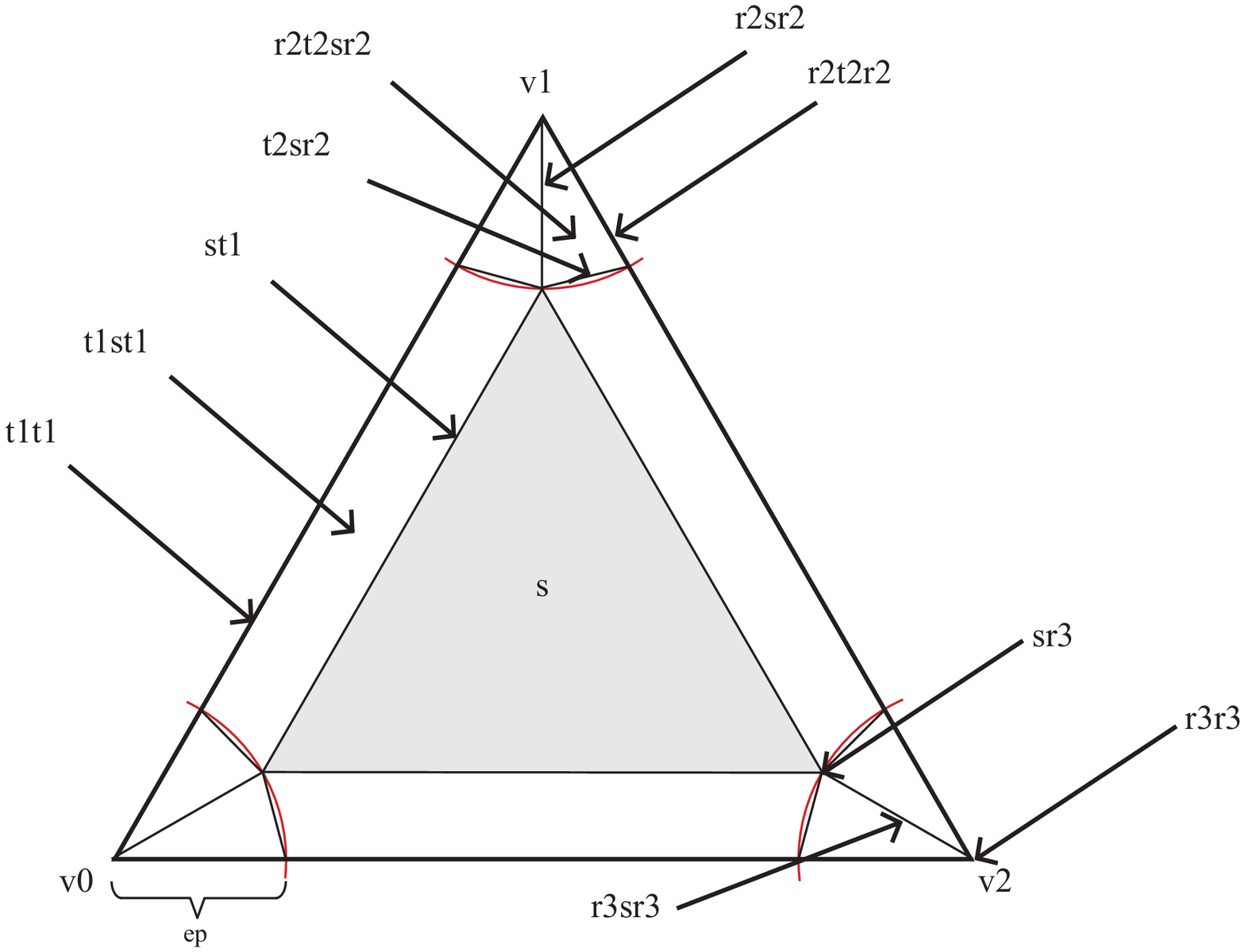}
}
\caption{The cellulation $X_\ep^\prime$ for a $2$-simplex.}
\label{Fig:fundamentalcellulation}
\end{center}
\end{figure}
\qed\end{ex}

\begin{rmk}\label{rmkcontrolep}
Note that for all $0<\ep<\comesh(X)$, $\Gamma_\ep$ is a p.l. isomorphism and that $\Gamma_0= \mathrm{pr}_1:X\times Sd\, X \to X$. Hence \[\Gamma_\delta \circ \Gamma_\ep^{-1}: X_\ep^\prime \to X_\delta^\prime\] is a p.l. isomorphism for all $0 < \ep,\delta < \comesh(X)$. Further, for $0<\ep<\comesh(X)$ the cellulation $X_\ep^\prime$ is homotopic to $X$ via the straight line homotopy
\begin{eqnarray*}
 h_{2,\ep}: Y\times I &\to& Y \\
 (y,t) &\mapsto& \Gamma_{\ep(1-t)}\Gamma_\ep^{-1}(y).
\end{eqnarray*}
This homotopy sends each vertex $\Gamma_\tau(v)$ to the point $v$ along a straight line of length precisely $\ep$. Convexity of the cells of $Y_\ep^\prime$ guarantees that all homotopy tracks are of length at most $\ep$. Hence $h_{2,\ep}$ has control $\ep$. 
\end{rmk}

\section{Proof of main theorem}\label{proof}
In this section we prove the main theorem which we restate for convenience.

\setcounter{thm1}{0}
\begin{thm1}
Let $f:X\to Y$ be a simplicial map of f.d. l.f. simplicial complexes equipped with their path metrics and let $j_Y:Y\times\R \to O(Y_+)$ be the coning map. Then the following are equivalent:
\begin{enumerate}[(i)]
 \item $f$ has contractible point inverses,
 \item $f:(X,f)\to (Y,\id_Y)$ is an $\ep$-controlled homotopy equivalence for all $\ep>0$,
 \item $f\times\id_\R: (X\times\R, j_Y(f\times\id_\R)) \to (Y\times\R,j_Y)$ is a bounded homotopy equivalence.
\end{enumerate}
\end{thm1}

To facilitate the proof of the main theorem we first require two propositions.

\begin{prop1}\label{propkaytaus}
Let $f:X\to Y$ be a simplicial map of l.f. f.d. simplicial complexes. Then 
\begin{enumerate}[(i)]
 \item for all simplices $\sigma \in Y$, there is a p.l. isomorphism $f^{-1}(\mathring{\sigma}) \cong f^{-1}(\widehat{\sigma})\times \mathring{\sigma},$
 \item $f$ has contractible point inverses if and only if $f^{-1}(\widehat{\sigma})$ is contractible for all $\sigma\in Y$.
\end{enumerate}
\end{prop1}

\begin{proof}
$(i)$: If $\mathring{\sigma}$ is not in the image of $f$ then the result holds as $f^{-1}(\mathring{\sigma})=f^{-1}(\widehat{\sigma})=\emptyset.$

Let $\sigma = w_0\ldots w_m$ be some simplex in $Y$. Suppose there is a $\tau\in X$ such that $f(\tau)=\sigma$. Let $f_\tau:= f|_\tau:\tau \to \sigma$. Since $\sigma$ is the join of its vertices we have that 
\[\tau = \moast_{i=0}^m{f_{\tau}^{-1}(w_i)}\]
with $f^{-1}_\tau(x) \cong \prod_{i=0}^{m}{f^{-1}_\tau(w_i)} \cong f_\tau^{-1}(\widehat{\sigma})$ for all $x\in \mathring{\sigma}$. Whence $f_\tau^{-1}(\mathring{\sigma}) \cong f_\tau^{-1}(\widehat{\sigma})\times \mathring{\sigma}.$

Suppose $\tau_0<\tau_1$ are such that $f(\tau_i)=\sigma$ for $i=0,1$. Then
\begin{displaymath}
 \xymatrix@R=3mm@C=3mm{
 f_{\tau_0}^{-1}(\mathring{\sigma}) \ar@{}[r]|-{\subset} \ar@{}[d]|-{\rotcong} & f_{\tau_1}^{-1}(\mathring{\sigma}) \ar@{}[d]|-{\rotcong} \\
 f_{\tau_0}^{-1}(\widehat{\sigma})\times\mathring{\sigma} \ar@{}[r]|-{\subset} & f_{\tau_1}^{-1}(\widehat{\sigma})\times\mathring{\sigma}.
 }
\end{displaymath}
Thus we can reconstruct $f^{-1}(\widehat{\sigma})$ from $\{ f^{-1}_\tau(\widehat{\sigma})| f(\tau)=\sigma\}$ as \[ f^{-1}(\widehat{\sigma}) = \bigcup_{\tau: f(\tau)=\sigma}{f^{-1}_\tau(\widehat{\sigma})}\] and consquently \[ f^{-1}(\mathring{\sigma}) = \bigcup_{\tau: f(\tau)=\sigma}{f^{-1}_\tau(\mathring{\sigma})} \cong \bigcup_{\tau: f(\tau)=\sigma}{f^{-1}_\tau(\widehat{\sigma})}\times \mathring{\sigma} = f^{-1}(\widehat{\sigma})\times\mathring{\sigma}.\]

$(ii)$: Clear from the fact that $f^{-1}(x)\cong f^{-1}(\widehat{\sigma})$ for $x\in \mathring{\sigma}$.
\end{proof}

This proposition tells us that for a simplicial map $f$ with contractible point inverses, the restriction over each simplex, $f|:f^{-1}(\mathring{\sigma})\to \mathring{\sigma}$, is a trivial fibre bundle with fibre $f^{-1}(\widehat{\sigma})\simeq *$. We will see that we can define a section over each simplex interior and the contractibility of each $f^{-1}(\widehat{\sigma})$ allows us to piece these local sections together by homotopies that are large in $X$ but can be made arbitrarily small in $Y$. This yields a global homotopy inverse $g_\ep$, for all $\ep>0$, that is an \textit{approximate section} in the sense that $f\circ g_\ep \simeq \id_Y$ via homotopy tracks of diameter $<\ep$. This approximate section can be used to approximately lift homotopies, hence we see that $f$ is an approximate fibration. 

\begin{prop1}\label{propdefret}
Let $f:X\to Y$ be a simplicial map of f.d. l.f. simplicial complexes. Then for all $\sigma\in Y$, $f^{-1}(\mathrm{st}(\sigma))$ p.l. deformation retracts onto $f^{-1}(\mathring{\sigma})$.
\end{prop1}
\begin{proof}
If $f^{-1}(\mathring{\sigma})$ is empty then so is $f^{-1}(\mathring{\rho})$ for all $\rho \geqslant \sigma$ and hence $f^{-1}(st(\sigma))$ is empty so the result holds vacuously.

Suppose instead that $f^{-1}(\mathring{\sigma})\cong f^{-1}(\widehat{\sigma})\times \mathring{\sigma}$ is non-empty. For every $\rho> \sigma$ let $\sigma_\rho^C\in Y$ be the unique simplex such that $\rho= \sigma * \sigma_\rho^C$. For all $\tau\in X$ with $f(\tau)=\rho$ let $f_\tau:= f|_\tau:\tau \to \rho$ so that $\tau = f^{-1}_\tau(\sigma)*f^{-1}_\tau(\sigma_\rho^C)$. Every $x\in \mathring{\tau}\cup f^{-1}_\tau(\mathring{\sigma})$ can be written uniquely as $x = (1-t)x_\sigma + tx_{\sigma_\tau^C}$, for $x_\sigma \in f^{-1}_\tau(\mathring{\sigma})$, $x_{\sigma_\rho^C}\in f^{-1}_\tau(\mathring{\sigma}_\rho^C)$, $t\in [0,1)$. Thus letting the $t$ parameter go to $0$ at unit speed and staying there thereafter defines a linear (strong) deformation retraction of $\mathring{\tau}\cup f^{-1}_\tau(\mathring{\sigma})$ onto $f^{-1}_\tau(\mathring{\sigma})$. The deformation retractions defined like this for different simplices surjecting onto $\rho$ agree on intersections and so glue to give a p.l. deformation retraction of $f^{-1}(\mathring{\rho}
\cup \mathring{\sigma})$ onto $f^{-1}(\mathring{\sigma})$. These glue together to give the desired deformation retraction of $f^{-1}(st(\sigma))$ onto $f^{-1}(\mathring{\sigma})$.
\end{proof}

\begin{proof}[Proof of Theorem \ref{maintopthm}]
$(i)\Rightarrow (iii)$:
Let $f:X\to Y$ be a simplicial map of f.d. l.f. simplicial complexes with contractible point inverses. Then $f$ is necessarily surjective as contractible point inverses are non-empty. We seek to define a one parameter family of homotopy inverses \[\{ g_\ep:Y\to X\}_{0<\ep<\comesh(Y)}\]and homotopies \[\{ h_{1,\ep}:\id_X \simeq g_\ep\circ f\}_{0<\ep<\comesh(Y)}, \quad \{ h_{2,\ep}:\id_Y \simeq f\circ g_\ep\}_{0<\ep<\comesh(Y)} \]parametrised by control. Given such families we obtain a bounded homotopy inverse $g$ to \[f\times \id_\R: (X\times\R,j_Y\circ (f\times\id_\R)) \to (Y\times\R,j_Y)\] defined by \[\begin{array}{rcl} g:Y\times\R &\to& X\times \R;\\ (y,t) &\mapsto& g_{\alpha(t)}(y)\end{array}\] and bounded homotopies  \[\begin{array}{rcl} h_1:\id_{X\times\R}\simeq g\circ(f\times\id_\R):X\times\R\times I &\to& X\times \R; \\ (x,t,s) &\mapsto& h_{1,{\alpha(t)}}(x,s),\\ h_2:\id_{Y\times\R}\simeq (f\times\id_\R)\circ g:Y\times\R\times I &\to& Y\times \R; \\ (y,t,s) &\mapsto& h_{2,{\alpha(t)}}(y,s),\end{array}
\] where $\alpha:\R \to (0,\comesh(Y)]$ is the function \[\alpha: t \mapsto \brcc{\comesh(Y),}{t\leqslant 1/\comesh(Y),}{1/t,}{t\geqslant 1/\comesh(Y).}\]

Give $Y$ the fundamental $\ep$-subdivision cellulation $Y_\ep^\prime$ as defined in Definition \ref{Defn:fundamentalcellulation}. We define $g_\ep$, $h_{1,\ep}$ and $h_{2,\ep}$ by induction. First, define a map $\gamma: \chi(Y) \to X$ by induction on the flag length of cells in $\chi(Y)$. Let \[\gamma_{\widehat{\sigma}\times \mathring{\sigma}}: \widehat{\sigma}\to f^{-1}(\widehat{\sigma})\] be any map, then define $\gamma$ on $\widehat{\sigma}\times \sigma$ as the closure of the map \[\gamma_{\widehat{\sigma}\times\mathring{\sigma}}\times \id_{\mathring{\sigma}}: \widehat{\sigma}\times \mathring{\sigma} \to f^{-1}(\widehat{\sigma})\times \mathring{\sigma} \cong f^{-1}(\mathring{\sigma}).\] 
Let $\Phi_{\tau,\sigma}:f^{-1}(\widehat{\sigma})\to f^{-1}(\widehat{\tau})$ denote the maps obtained in the closure of $\gamma_{\widehat{\sigma}\times\mathring{\sigma}}$ for $\tau<\sigma$ such that \[\gamma_{\widehat{\sigma}\times \mathring{\tau}} = (\Phi_{\tau,\sigma}\circ \gamma_{\widehat{\sigma}\times \mathring{\sigma}})\times \id_{\mathring{\tau}} : \widehat{\sigma}\times \mathring{\tau} \to f^{-1}(\widehat{\tau})\times \mathring{\tau} \cong f^{-1}(\mathring{\tau}).\]

Now suppose that we have continuously defined $\gamma$ on all cells of $\chi(Y)$ of flag length at most $n$ and that the map takes the form \[\gamma_{\widehat{\sigma}_0\ldots \widehat{\sigma}_i \times \mathring{\sigma}_0}\times \id_{\mathring{\sigma}_0}: \widehat{\sigma}_0\ldots \widehat{\sigma}_i \times \mathring{\sigma}_0 \to f^{-1}(\mathring{\sigma}_0)\times \mathring{\sigma}_0\cong f^{-1}(\mathring{\sigma}_0)\] on each cell for $i\leqslant n$ for some maps \[\gamma_{\widehat{\sigma}_0\ldots \widehat{\sigma}_i \times \mathring{\sigma}_0}: \widehat{\sigma}_0\ldots \widehat{\sigma}_i \to f^{-1}(\mathring{\sigma}_0).\] These maps define a map \[\gamma_{\partial(\widehat{\sigma}_0\ldots \widehat{\sigma}_{n+1}) \times \mathring{\sigma}_0}: \partial(\widehat{\sigma}_0\ldots \widehat{\sigma}_{n+1}) \to f^{-1}(\widehat{\sigma}_0)\] which extends to a map \[\gamma_{\widehat{\sigma}_0\ldots \widehat{\sigma}_{n+1}\times \mathring{\sigma}_0}: \widehat{\sigma}_0\ldots \widehat{\sigma}_{n+1} \to f^{-1}(\widehat{\sigma}_0)\] by the contractibility of $f^{-1}(\widehat{\sigma}_0)$. Define $\gamma$ on the cell $\widehat{\sigma}_0\ldots \widehat{\sigma}_{n+1}\times \sigma_0$ as the closure of the map \[ \gamma_{\widehat{\sigma}_0\ldots \widehat{\sigma}_{n+1}\times \mathring{\sigma}_0}\times \id_{\mathring{\sigma}_0}: \gamma_{\widehat{\sigma}_0\ldots \widehat{\sigma}_{n+1}\times \mathring{\sigma}_0} \to f^{-1}(\widehat{\sigma}_0)\times \mathring{\sigma}_0\cong f^{-1}(\mathring{\sigma}_0).\]
By induction this defines the map $\gamma$.

For all $0<\ep<\comesh(Y)$, set \[g_\ep:= \gamma\circ \Gamma_\ep^{-1}: Y \to \chi(Y) \to X.\] We claim that $\{g_\ep\}_{0<\ep<\comesh(Y)}$ is a one parameter family of homotopy inverses to $f$ parametrised by control.

Consider first the composition $f\circ g_\ep$. \[f\circ \gamma = pr_2\circ (\gamma_{\widehat{\sigma}_0\ldots \widehat{\sigma}_{n}\times \mathring{\sigma}_0}\times \id_{\mathring{\sigma}_0}) = pr_2: \widehat{\sigma}_0\ldots \widehat{\sigma}_{n}\times \mathring{\sigma}_0 \to f^{-1}(\widehat{\sigma}_0)\times \mathring{\sigma}_0 \to \mathring{\sigma}_0.\] Hence $f\circ \gamma: \chi(Y)\subset Sd\, Y \times Y \to Y$ is just projection onto $Y$, i.e. the map $\Gamma_0 = \lim_{\ep \to 0}{\Gamma_\ep}$. Thus $f\circ g_\ep = (f\circ \gamma)\circ \Gamma^{-1}_\ep = \Gamma_0\circ \Gamma^{-1}_\ep$. Choosing $h_{2,\ep}$ precisely as in Remark \ref{rmkcontrolep} we have $h_{2,\ep}: \id_X = \Gamma_\ep\circ \Gamma^{-1}_\ep \simeq \Gamma_0\circ \Gamma^{-1}_\ep= f\circ g_\ep$ is an $\ep$-controlled homotopy and in fact \[\{h_{2,\ep}\}_{0<\ep<\comesh(Y)}\] is a one parameter of homotopies parametrised by control.

Now consider the other composition: $g_\ep\circ f = \gamma\circ \Gamma_\ep^{-1} \circ f$. Define a homotopy $h_{1,\ep}^\prime: X\times I \to X$ by \[h_{1,\ep}^\prime = \id_{f^{-1}(\widehat{\sigma})}\times h_{2,\ep}: f^{-1}(\widehat{\sigma})\times \mathring{\sigma}\times [0,1) \to f^{-1}(\mathring{\sigma})\] with $h_{1,\ep}^\prime(-,1):= \lim_{t\to 1}{h_{1,\ep}^\prime(-,t)}.$ This homotopy is sent by $f$ to $h_{2,\ep}$: \[ f(h_{1,\ep}^\prime(x,t)) = h_{2,\ep}(f(x),t),\quad \forall (x,t)\in X\times I.\] Hence $h_{1,\ep}^\prime$ has control $\ep$.

We now seek a homotopy $h_{1,\ep}^{\prime\prime}:h_{1,\ep}^\prime(-,1)\simeq g_\ep\circ f$ with zero control. Looking at $f^{-1}(\Gamma_\ep(\rho\times \mathring{\sigma}_0))$ for $\rho = \widehat{\sigma}_0\ldots \widehat{\sigma}_{n}$ observe that $h_{1,\ep}^\prime(-,1)$ is the closure of the map \[\Phi_{\sigma_0,\sigma_n} \times h_{2,\ep}(-,1) =  \Phi_{\sigma_0,\sigma_n} \times (\Gamma_0\circ \Gamma_\ep^{-1}): f^{-1}(\widehat{\sigma}_n) \times \Gamma_\ep(\mathring{\rho}\times \mathring{\sigma}_0) \to f^{-1}(\widehat{\sigma}_0)\times \mathring{\sigma}_0,\] whereas $g_\ep\circ f$ is the closure of the map \[ (\gamma_{\widehat{\sigma}_0\ldots \widehat{\sigma}_{n}\times \mathring{\sigma}_0} \times \id_{\mathring{\sigma}_0})\circ \Gamma_\ep^{-1} \circ pr_2: f^{-1}(\widehat{\sigma}_n) \times \Gamma_\ep(\mathring{\rho}\times \mathring{\sigma}_0) \to f^{-1}(\widehat{\sigma}_0)\times \mathring{\sigma}_0.\] The component of this map from $\Gamma_\ep(\mathring{\rho})$ to $\mathring{\sigma}_0$ is $\Gamma_0\Gamma_\ep^{-1}$ and so agrees 
with the component of $h_{1,\ep}^\prime(-,1)$ to $\mathring{\sigma}_0$. We now find inductively a homotopy $h_{1,\ep}^{\prime\prime}:h_{1,\ep}^\prime(-,1)\simeq g_\ep\circ f$ which only moves things in the fibre direction and hence has $0$ control. This is achieved precisely as before using the contractibility of the fibres. The concatenation $h_{1,\ep}:= h_{1,\ep}^{\prime\prime}*h_{1,\ep}^{\prime}$ is an $\ep$-controlled homotopy $\id_X\simeq g_\ep\circ f$. As we use the same homotopies in the fibre direction for all $0<\ep<\comesh(Y)$ this gives a one parameter family $\{ h_{1,\ep}:\id_Y\simeq g_\ep\circ f\}_{0<\ep<\comesh(Y)}$ parametrised by control as required.

Note also that $f|:f^{-1}(\tau)\to \tau$ is a homotopy equivalence for all $\tau\in Y$ by restricting $g_\ep, h_{1,\ep}$ and $h_{2,\ep}$. We call such a homotopy equivalence a \textit{$Y$-triangular homotopy equivalence}. It is an open conjecture that $f:X\to Y$ is homotopic to a $Y$-triangular homotopy equivalence if and only if $f$ is homotopic to an $\ep$-controlled homotopy equivalence for all $\ep>0$. $Y$-triangular homotopy equivalences are discussed in \cite{shortpaper}.

$(iii)\Rightarrow (ii)$:
Let $f\times\id$ have homotopy inverse $g$ and homotopies $h_1:\id_{X\times\R}\simeq g\circ (f\times\id_\R) $ and $h_2:\id_{Y\times\R} \simeq (f\times\id_\R)\circ g$ all with bound at most $B<\infty$. Let $p_t:\R\to \{t\}$ be projection onto $t\in \R$. 

Let $g_t:= (\id_X\times p_t)\circ g|_{Y\times\{t\}}: Y\times\{t\} \to X\times\R \to X\times\{t\}.$ This is a homotopy inverse to $f\times\id_{\{t\}}:X\times\{t\}\to Y\times\{t\}$ with homotopies
\begin{align}
 (\id_X\times p_t)\circ h_1|_{X\times\{t\}}: (\id_X\times p_t)\circ \id_{X\times\{t\}} = \id_{X\times\{t\}} &\simeq (\id_X\times p_t)\circ (g\circ (f\times\id_\R))|_{X\times\{t\}} \notag \\ &= (\id_X\times p_t)\circ g|_{Y\times\{t\}} \circ (f\times\id_{\{t\}}) \notag \\ &= g_t\circ (f\times\id_{\{t\}})\notag
\end{align}
and
\begin{align}
 (\id_Y\times p_t)\circ h_2|_{Y\times\{t\}}: (\id_Y\times p_t)\circ \id_{Y\times\{t\}} = \id_{Y\times\{t\}} &\simeq (\id_Y\times p_t)\circ ((f\times\id_\R)\circ g)|_{Y\times\{t\}} \notag \\ &= (\id_Y\times p_t)\circ (f\times\id_\R) \circ g|_{Y\times\{t\}}  \notag \\ &= (f\times\id_{\{t\}})\circ (\id_X\times p_t)\circ g|_{Y\times\{t\}} \notag \\ &= (f\times\id_{\{t\}})\circ g_t. \notag
\end{align}

These homotopies have bound approximately $B$ measured in $Y\times\{t\} \subset O(Y_+)$. The slice $Y\times\{t\}$ has a metric $t$ times bigger than $Y=Y\times\{1\}$, so measuring this in $Y$ gives a homotopy equivalence $f:X\to Y$ with control proportional to $\dfrac{B}{t}$ as required.

$(ii)\Rightarrow(i)$:
First note that a simplicial map $f$ that is an $\ep$-controlled homotopy equivalence for all $\ep>0$ must be surjective. Suppose it is not, then there is a $y\in Y\backslash \mathrm{im}(f)$. Since $f$ is simplicial $\mathring{\sigma}\subset Y\backslash \mathrm{im}(f)$ where $\sigma$ is the unique simplex of $Y$ with $y\in \mathring{\sigma}$. Again since $f$ is simplicial, if $\tau\geqslant \sigma$ we must have $\mathring{\tau} \subset Y\backslash \mathrm{im}(f)$ as well. Thus \[st(\sigma) = \bigcup_{\tau\geqslant\sigma}{\mathring{\tau}} \;\subset\; Y\backslash \mathrm{im}(f).\] In particular the open star $st(\sigma)$ is an open neighbourhood of $y$ in $Y\backslash \mathrm{im}(f)$ so we may find a ball $B_{\ep^\prime}(y)\subset Y\backslash \mathrm{im}(f)$. Thus $f$ cannot be an $\ep$-controlled homotopy equivalence for $\ep<\ep^\prime$ as the homotopy tracks for the point $y$ must travel a distance of at least $\ep^\prime$. This is a contradiction and so $f$ is surjective.

Each point $y\in Y$ is contained in a unique simplex interior and hence in that simplex's open star: $\mathring{\sigma}\subset\mathrm{st}(\sigma)$. Since the star is open there is an $\ep^\prime$ such that $B_{\ep^\prime}(y)\subset\mathrm{st}(\sigma)$. By hypothesis we can find an $\ep^\prime$-controlled homotopy inverse, $g_{\ep^\prime}$, to $f$. Thus $f^{-1}(y)$ is homotopic to $g_{\ep^\prime}(y)$ within $f^{-1}(\mathrm{st}(\sigma))$. By Proposition \ref{propdefret}, $f^{-1}(\mathrm{st}(\sigma))$ deformation retracts onto $f^{-1}(\mathring{\sigma})$. By Proposition \ref{propkaytaus} this is p.l. isomorphic to $f^{-1}(\widehat{\sigma})\times \mathring{\sigma}$ which in turn deformation retracts onto $f^{-1}(\widehat{\sigma})\times \{y\} = f^{-1}(y)$. Applying these two deformation retractions to the homotopy $f^{-1}(y)\simeq g_{\ep^\prime}(y)$ gives a contraction of $f^{-1}(y)$. Hence $f$ has contractible point inverses.
% OLD PROOF
% By Lemma \ref{fsurjects} $f$ must necessarily be a surjective map. By Lemma \ref{surjectivesimplicialmapssimplextosimplex}, for all $\sigma\in Y$ \[f^{-1}(\mathring{\sigma})\cong \mathring{\sigma}\times K(\sigma)\] for some simplicial complex $K(\sigma)$. Pick $x\in \mathring{\sigma}$, then there exists an $\ep$ such that $x\in B_{\ep}(x)\subset \mathring{\sigma}$, whence $f^{-1}(B_\ep(x))\cong B_\ep(x)\times K(\sigma)$. Since $f$ is an $\ep$-controlled homotopy equivalence for all $\ep$ we can find an $\ep/2$-controlled homotopy inverse $g$ together with homotopies such that $f^{-1}(x) \simeq g(x)=*$ in $f^{-1}(B_\ep(x))$. Projecting this homotopy to $K(\sigma)$ gives a contraction of $K(\sigma)$ to $\pr_2(g(x))=*$. Thus $K(\sigma)\simeq *$ for all $\sigma\in Y$. Applying Lemma \ref{Ksigmacontractible} $(iii)$ we deduce that $f$ is a contractible map.
\end{proof}

We conclude with an example illustrating the construction in the proof of $(i)\Rightarrow (iii)$.
\begin{ex}
Let $\underline{0}=(0,0,0)$, $e_1=(1,0,0)$, $e_2 = (0,1,0)$ and $e_3=(0,0,1)$ be points in $\R^3$. Define $Y$ to be the simplicial complex with the following $2$-simplices: $\sigma_1:=\underline{0}*e_1*(e_1+e_2)$ and $\sigma_2:=\underline{0}*e_2*(e_1+e_2)$. Define $X$ to be the simplicial complex with the following $2$-simplices: $\tau_1:=\underline{0}*e_1*(e_1+e_2)$, $\tau_2:=e_3*(e_2+e_3)*(e_1+e_2+e_3)$, $\tau_3:=\underline{0}*e_3*(e_1+e_2+e_3)$ and $\tau_4:=\underline{0}*(e_1+e_2)*(e_1+e_2+e_3)$. The projection map $f:X\to Y;(x,y,z)\mapsto (x,y,0)$ is simplicial and has contractible point inverses. Give $Y$ the cellulation $Y_\ep^\prime$ for some small $\ep>0$ as pictured in \Figref{Fig:Workedex1}. 

\begin{figure}[tbh!]
\begin{center}
{
 \psfrag{x}[]{$x$}
 \psfrag{y}[]{$y$}
% \psfrag{t=1}[]{$t=1$}
\includegraphics[width=6cm]{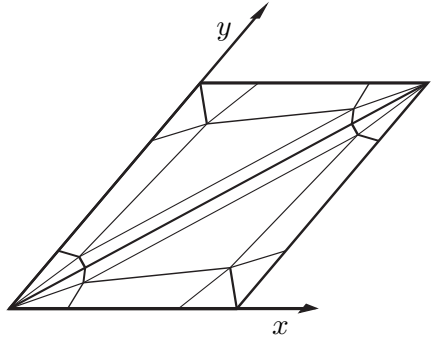}
}
\caption{$\ep$-subdivision cellulations.}
\label{Fig:Workedex1}
\end{center}
\end{figure}

We define $g_\ep$ as in the proof by first defining maps $\gamma_{\widehat{\rho}\times \mathring{\rho}}: \widehat{\rho}\to f^{-1}(\widehat{\rho})$ for all $\rho\in Y$. We define 
\[\gamma_{\widehat{\rho}\times \mathring{\rho}}= \left\{ \begin{array}{cc} 0, & \mathring {\rho}\subset \sigma_1\backslash \sigma_2, \\ 1/2, & \mathring {\rho}\subset \sigma_1\cap \sigma_2, \\ 1, & \mathring {\rho}\subset \sigma_2\backslash \sigma_1.\end{array}\right.\] Then, for all $\rho\in\sigma_1\cap \sigma_2$ we choose the maps 
\[ \gamma_{\widehat{\rho}\widehat{\sigma}_i\times \mathring{\rho}}:\widehat{\rho}\widehat{\sigma}_i \to f^{-1}(\widehat{\rho})\]
for $i=1,2$ as follows:
\begin{align}
 \gamma_{\widehat{\rho}\widehat{\sigma}_1\times \mathring{\rho}}(t_0,t_1) &= \frac{1}{2}t_0, \notag \\
 \gamma_{\widehat{\rho}\widehat{\sigma}_2\times \mathring{\rho}}(t_0,t_1) &= \frac{1}{2}t_0+ t_1 \notag
\end{align}
where $(t_0,t_1)$ are barycentric coordinates.

Finally for $v$ either vertex of $\sigma_1\cap \sigma_2$ and $\rho$ the $1$-simplex of $\sigma_1\cap \sigma_2$ we define the maps 
\[ \gamma_{\widehat{v}\widehat{\rho}\widehat{\sigma}_i\times \mathring{v}}:\widehat{v}\widehat{\rho}\widehat{\sigma}_i \to f^{-1}(\widehat{v})\]
for $i=1,2$ as follows:
\begin{align}
 \gamma_{\widehat{v}\widehat{\rho}\widehat{\sigma}_1\times \mathring{v}}(t_0,t_1,t_2) &= \frac{1}{2}t_0 + \frac{1}{2}t_1, \notag \\
 \gamma_{\widehat{v}\widehat{\rho}\widehat{\sigma}_2\times \mathring{v}}(t_0,t_1,t_2) &= \frac{1}{2}t_0 + \frac{1}{2}t_1 + t_2. \notag
\end{align}
The resulting map $g_\ep$ is illustrated in \Figref{Fig:Workedex2} where it is exaggerated to show where each cell of $Y_\ep^\prime$ is sent.
\begin{figure}[tbh!]
\begin{center}
{
 \psfrag{x}[]{$x$}
 \psfrag{y}[]{$y$}
 \psfrag{z}[]{$z$}
 \psfrag{0.5}[]{$0.5$}
 \psfrag{1}[]{$1\;$}
 \psfrag{0}[]{$0\;$}
\includegraphics[width=7.5cm]{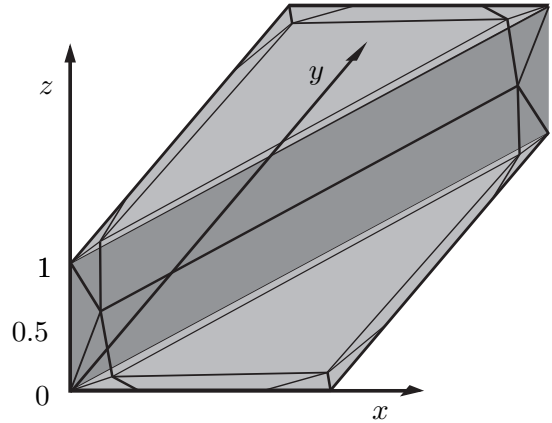}
}
\caption{Constructing the map $g_\ep$.}
\label{Fig:Workedex2}
\end{center}
\end{figure}
\end{ex}

\bibliographystyle{hep}
\bibliography{spirosbib}{}
\end{document}